\documentclass[12pt,leqno,a4paper]{article}
\usepackage{amsmath}
\usepackage{amsfonts}
\usepackage{amsthm}
\usepackage{amssymb}

\usepackage{xspace}
\usepackage{euscript}
\usepackage{graphicx}
\usepackage{amscd}
\usepackage{tabularx}

\usepackage{enumerate}
 \usepackage{epsfig} 
 \usepackage{graphics} 
\theoremstyle{plain}
\newtheorem{theo}{Theorem}[section]
\newtheorem{prop}[theo]{Proposition}
\theoremstyle{definition}
\newtheorem{definition}[theo]{Definition}
\theoremstyle{remark}
\newtheorem{rem}[theo]{Remark}
\numberwithin{equation}{section}

\newcommand{\C}{\mathbb{C}}

\newcommand{\K}{\mathbb{K}}
\newcommand{\R}{\mathbb{R}}
\newcommand{\N}{\mathbb{N}}
\newcommand{\M}{\mathbb{M}}
\newcommand{\divrg}{\textrm{div}\,}

\title{Optimal identification of cavities in the Generalized Plane Stress problem in linear elasticity}
\author{Antonino Morassi\thanks{Dipartimento Politecnico di Ingegneria e Architettura,
Universit\`a degli Studi di Udine, via Cotonificio 114, 33100
Udine, Italy. E-mail: \textsf{antonino.morassi@uniud.it}}, \  Edi
Rosset\thanks{Dipartimento di Matematica e Geoscienze,
Universit\`a degli Studi di Trieste, via Valerio 12/1, 34127
Trieste, Italy. E-mail: \textsf{rossedi@units.it}} \ and Sergio
Vessella\thanks{Dipartimento di Matematica e Informatica ``Ulisse
Dini'', Universit\`a degli Studi di Firenze, Via Morgagni 67/a,
50134 Firenze, Italy. E-mail: \textsf{sergio.vessella@unifi.it}}}

\begin{document}

\maketitle

\begin{abstract}

For the Generalized Plane Stress (GPS) problem in linear
elasticity, we obtain an optimal stability estimate of logarithmic
type for the inverse problem of determining smooth cavities inside
a thin isotropic cylinder {}from a single boundary measurement of
traction and displacement. The result is obtained by reformulating
the GPS problem as a Kirchhoff-Love plate-like problem in terms of
the Airy's function, and by using the strong unique continuation
at the boundary for a Kirchhoff-Love plate operator under
homogeneous Dirichlet conditions, which has been recently obtained
in \cite{l:arv}.

\medskip


\medskip

\end{abstract}

\centerline{}

\section{Introduction}
\label{Introduction}

In this paper we consider the inverse problem of detecting
cavities inside a thin isotropic elastic plate $\Omega \times
\left (  -\frac{h}{2}, \frac{h}{2} \right )$, where the middle
plane $\Omega$ is a bounded domain in $\R^2$ and $h$ is the
constant thickness, subject to a single experiment consisting in
applying in-plane boundary loads and measuring the induced
displacement at the boundary. Practical applications concern the
use of non-destructive techniques for the identification of
possible defects, such as cavities, inside the plate.

The static equilibrium of the plate is described in terms of the
classical Generalized Plane Stress (GPS) problem, which allows to
reformulate the original three dimensional problem in a two
dimensional setting \cite{Sokolnikoff}. More precisely, denoting
by $D \times \left (  -\frac{h}{2}, \frac{h}{2} \right )$ the
cavity, with $D$ a possibly disconnected subset of $\Omega$, the
in-plane displacement field $a=a_1 e_1+a_2 e_2$, solution to the
GPS problem, satisfies the following two-dimensional Neumann
boundary value problem $(\alpha, \beta =1,2)$
\begin{center}
\( {\displaystyle \left\{
\begin{array}{lr}
     N_{\alpha\beta,\beta}=0,
      & \mathrm{in}\ \Omega \setminus \overline {D},
        \vspace{0.25em}\\
      N_{\alpha\beta}n_\beta=\widehat{N}_\alpha, & \mathrm{on}\ \partial \Omega,
          \vspace{0.25em}\\
      N_{\alpha\beta}n_\beta=0, & \mathrm{on}\ \partial D,
        \vspace{0.25em}\\
      N_{\alpha\beta}=\frac{Eh}{1-\nu^2}\left((1-\nu)\epsilon_{\alpha\beta}+\nu(\epsilon_{\gamma\gamma})\delta_{\alpha\beta}\right), &\mathrm{in}\ \Omega \setminus \overline {D},
          \vspace{0.25em}\\
        \epsilon_{\alpha\beta}=\frac{1}{2}(a_{\alpha,\beta}+a_{\beta,\alpha}), &\mathrm{in}\ \Omega \setminus \overline {D}.
          \vspace{0.25em}\\
\end{array}
\right. } \) \vskip -8.9em
\begin{eqnarray}
& & \label{eq:I1}\\
& & \label{eq:I2}\\
& & \label{eq:I3}\\
& & \label{eq:I4}\\
& & \label{eq:I5}
\end{eqnarray}
\end{center}
Here, $\widehat{N}=\widehat{N}_1e_1+\widehat{N}_2e_2$ is the
in-plane load field applied to $\partial\Omega$ satisfying the
compatibility condition
\begin{equation}
  \label{eq:I6}
  \int_{\partial\Omega}\widehat{N}\cdot r =0, \quad \hbox{for every } r\in {\mathcal R_2},
\end{equation}
where ${\mathcal R_2}$ is the linear space of infinitesimal
two-dimensional rigid displacements. Here, $E=E(x)$ and
$\nu=\nu(x)$ are the Young's modulus and the Poisson's coefficient
of the material, respectively. Under suitable strong convexity
assumptions on the elastic tensor of the material (see Section
\ref{sec:IP} for details), and assuming
 $\widehat{N}\in
H^{-\frac{1}{2}}(\partial\Omega,\R^2)$, problem
\eqref{eq:I1}--\eqref{eq:I6} admits a unique solution $a \in
H^1(\Omega\setminus \overline{D},\R^2)$ satisfying the
normalization conditions
\begin{equation}
  \label{eq:I7}
\int_{\Omega \setminus \overline {D}}a =0, \quad \int_{\Omega
\setminus \overline {D}}(\nabla a - \nabla^Ta) =0,
\end{equation}
and such that $\|a\|_{H^1(\Omega \setminus \overline {D})}\leq
C\|\widehat{N}\|_{H^{-\frac{1}{2}}(\partial\Omega,\R^2)}$.

In this work we face the inverse problem of determining the cavity
$D$ {}from a single pair of Cauchy data $\{a,\widehat{N}\}$ given
on $\partial \Omega$. More precisely, we are interested to obtain
quantitative stability estimates, which are useful to control the
effect that possible errors on the measurements have on the
results of reconstruction procedures. The arbitrariness of the
normalization conditions \eqref{eq:I7}, which are related to the
non-uniqueness of the solution to the direct problem
\eqref{eq:I1}--\eqref{eq:I6}, leads to the following formulation
of the stability issue: given two solutions $a^{(i)} \in
H^1(\Omega, \R^2)$, $i=1,2$, to the direct problem
\eqref{eq:I1}--\eqref{eq:I6} with $D=D_i$, satisfying, for some
$\varepsilon >0$,
\begin{equation}
    \label{eq:I8}
\min_{r \in \mathcal R_2} \|a^{(1)} - a^{(2)} -r \|_{L^2(\Sigma,
\R^2)}\leq \varepsilon,
\end{equation}
to control the Hausdorff distance $d_{\cal H}(
\overline{D_1},\overline{D_2} )$ in terms of $\varepsilon$ when
$\varepsilon$ goes to zero, where $\Sigma$ is an open subset of
$\partial \Omega$.

Assuming $D \in C^{6,\alpha}$, $0 < \alpha \leq 1$, we prove
\begin{equation}
    \label{eq:I9}
d_{\cal H}( \overline{D_1},\overline{D_2} ) \leq C | \log
\varepsilon |^{-\eta},
\end{equation}
where $C>0$ and $\eta>0$ are constants only depending on the a
priori data. We refer to Theorem \ref{theo:Main} for a precise
statement. Let us notice that, in view of the counterexamples
obtained in the simpler context of electrical conductivity (see,
for instance, \cite{l:a2}, \cite{Mandache}, \cite{l:dr}), we can
infer the optimality of the stability estimate \eqref{eq:I9}.

The general scheme of our proof is inspired to the seminal paper
\cite{A-B-R-V}, which established the first optimal logarithmic
estimate for the determination of unknown boundaries in
electrostatics. The key tool in \cite{A-B-R-V} was, among others,
the polynomial vanishing rate for solutions to the second order
elliptic equation of electrostatics, satisfying either homogeneous
Dirichlet or homogeneous Neumann boundary conditions, ensured by a
doubling inequality at the boundary established in \cite{l:ae}.
Aiming at obtaining a strong unique continuation property at the
boundary (SUCB) for solutions to the GPS elliptic system, in this
paper we have exploited the two dimensional character of the
problem \eqref{eq:I1}--\eqref{eq:I6} by using the classical Airy's
transformation, which (locally) reduces the GPS system with
homogeneous Neumann boundary conditions to a scalar fourth order
Kirchhoff-Love plate's equation under homogeneous Dirichlet
boundary conditions. This reformulation allows us to use the
finite vanishing rate at the boundary for homogeneous Dirichlet
boundary conditions recently obtained in \cite{l:arv} in the form
of a three spheres inequality at the boundary with optimal
exponent, and in \cite{MRV-Doubl-LeMat2019} in the form of a
doubling inequality at the boundary.

It is worth noticing that the present approach, here applied to
the GPS problem, allows also to cover the analogous inverse
problem of detecting cavities in a two-dimensional elastic body
made by inhomogeneous Lam\'{e} material, thus improving the
$\log-\log$ stability result previously obtained in
\cite{l:mr2004}. An optimal log-type estimate in dimension three
remains a challenging open problem. Let us mention that the Airy's
transformation has been used in \cite{L-U-W} to prove global
identifiability of the viscosity in an incompressible fluid
governed by the Stokes and the Navier-Stokes equations in the
plane by using boundary measurements.

The paper is organized as follows. Notation is presented in
Section \ref{sec:notation}. Section \ref{sec:IP} contains the
formulation of the inverse problem and the statement of our
stability result. The Airy's transformation is illustrated in
Section \ref{sec:Airy}. The proof of the main result, given in
Section \ref{sec:proof}, is based on a series of auxiliary
propositions concerning Lipschitz propagation of smallness
(Proposition \ref{prop:LPS}), finite vanishing rate in the
interior (Proposition \ref{prop:FVRI}), finite vanishing rate at
the boundary (Proposition \ref{prop:FVRB}), stability estimate
{}from Cauchy data (Proposition \ref{prop:4.8.1}). Finally, for
the sake of completeness, in Section \ref{sec:GPS} we recall a
derivation of the GPS problem {}from the corresponding three
dimensional elasticity problem for a thin plate subject to
in-plane boundary loads.

\section{Notation}
   \label{sec:notation}

Let $P=(x_1(P), x_2(P))$ be a point of $\R^2$. We shall denote by
$B_r(P)$ the disk in $\R^2$ of radius $r$ and center $P$ and by
$R_{a,b}(P)$ the rectangle of center $P$ and sides parallel to the
coordinate axes, of length $2a$ and $2b$, namely
$R_{a,b}(P)=\{x=(x_1,x_2)\ |\ |x_1-x_1(P)|<a,\ |x_2-x_2(P)|<b \}$.
\begin{definition}
  \label{def:reg_bordo} (${C}^{k,\alpha}$ regularity)
Let $\Omega$ be a bounded domain in ${\R}^{2}$. Given $k,\alpha$,
with $k\in\N$, $0<\alpha\leq 1$, we say that a portion $S$ of
$\partial \Omega$ is of \textit{class ${C}^{k,\alpha}$ with
constants $r_{0}$, $M_{0}>0$}, if, for any $P \in S$, there exists
a rigid transformation of coordinates under which we have $P=0$
and
\begin{equation*}
  \Omega \cap R_{r_0,2M_0r_0}=\{x \in R_{r_0,2M_0r_0} \quad | \quad
x_{2}>g(x_1)
  \},
\end{equation*}
where $g$ is a ${C}^{k,\alpha}$ function on $[-r_0,r_0]$
satisfying
\begin{equation*}
g(0)=g'(0)=0,
\end{equation*}
\begin{equation*}
\|g\|_{{C}^{k,\alpha}([-r_0,r_0])} \leq M_0r_0,
\end{equation*}
where
\begin{equation*}
\|g\|_{{C}^{k,\alpha}([-r_0,r_0])} = \sum_{i=0}^k
r_0^i\sup_{[-r_0,r_0]}|g^{(i)}|+r_0^{k+\alpha}|g|_{k,\alpha},
\end{equation*}

\begin{equation*}
|g|_{k,\alpha}= \sup_ {\overset{\scriptstyle t,s\in
[-r_0,r_0]}{\scriptstyle t\neq s}}\left\{\frac{|g^{(k)}(t) -
g^{(k)}(s)|}{|t-s|^\alpha}\right\}.
\end{equation*}

\end{definition}

We use the convention to normalize all norms in such a way that
their terms are dimensionally homogeneous and coincide with the
standard definition when the dimensional parameter equals one. For
instance,
\begin{equation*}
    \|f\|_{H^1(\Omega)}=r_0^{-1} \left ( \int_\Omega f^2
    +r_0^2\int_\Omega|\nabla f|^2 \right )^{\frac{1}{2}},
\end{equation*}
and so on for boundary and trace norms.

Given a bounded domain $\Omega$ in $\R^2$ such that $\partial
\Omega$ is of class $C^{k,\alpha}$, with $k\geq 1$, we consider as
positive the orientation of the boundary induced by the outer unit
normal $n$ in the following sense. Given a point
$P\in\partial\Omega$, let us denote by $\tau=\tau(P)$ the unit
tangent at the boundary in $P$ obtained by applying to $n$ a
counterclockwise rotation of angle $\frac{\pi}{2}$, that is
\begin{equation}
    \label{eq:2.tangent}
        \tau=e_3 \times n,
\end{equation}
where $\times$ denotes the vector product in $\R^3$ and $\{e_1,
e_2, e_3\}$ is the canonical basis in $\R^3$.

Given any connected component $\cal C$ of $\partial \Omega$ and
fixed a point $P_0\in\cal C$, let us define as positive the
orientation of $\cal C$ associated to an arclength
parameterization $\psi(s)=(x_1(s), x_2(s))$, $s \in [0, l(\cal
C)]$, such that $\psi(0)=P_0$ and $\psi'(s)=\tau(\psi(s))$. Here
$l(\cal C)$ denotes the length of $\cal C$.

Throughout the paper, we denote by $w,_\alpha$, $\alpha=1,2$,
$w,_s$, and $w,_n$ the derivatives of a function $w$ with respect
to the $x_\alpha$ variable, to the arclength $s$ and to the normal
direction $n$, respectively, and similarly for higher order
derivatives.

We denote by $\M^{n}$ the space of $n \times n$ real valued
matrices and by ${\cal L} (X, Y)$ the space of bounded linear
operators between Banach spaces $X$ and $Y$.

Given $A$, $B \in \M^{n}$ and $\K\in{\cal L} ({\M}^{n},
{\M}^{n})$, we use the following notation:
\begin{equation}
  \label{eq:2.notation_1}
  ({\K}A)_{ij} = \sum_{k,l=1}^{n} K_{ijkl}A_{kl},
\end{equation}
\begin{equation}
  \label{eq:2.notation_2}
  A \cdot B = \sum_{i,j=1}^{n} A_{ij}B_{ij},
\end{equation}
\begin{equation}
  \label{eq:2notation_3}
  |A|= (A \cdot A)^{\frac {1} {2}},
\end{equation}
\begin{equation}
  \label{eq:2.notation_4}
  \widehat {A} = \frac {1} {2} (A + A^{T}).
\end{equation}
We denote by $I_n$ the $n\times n$ identity matrix, and by
$\hbox{tr}(A)$ the trace of $A$.

When $n=2$, we replace the Latin indexes with Greek ones.

The linear space of the infinitesimal rigid displacements, for
$n=2, \ 3$, is defined as
\begin{equation}
  \label{eq:def_rig_displ-1}
    {\cal R}_n = \left \{
    r(x) = c + Wx, \ c \in \R^n, \ W \in \M^n, \ W+W^T=0
    \right \}.
\end{equation}

\section{Inverse problem and main result}
   \label{sec:IP}

\noindent {\it i) A priori information on the geometry.}

Let $\Omega$ be a bounded domain in $\R^2$ and let us assume that
the cavity $D$ is an open subset compactly contained in $\Omega$,
such that
\begin{equation}
  \label{eq:4.1.1}
  \Omega\setminus D \ \rm{is\ connected}.
\end{equation}
Moreover, let us assume that, given positive numbers $r_0$, $M_0$, $M_1$, with
$M_0\geq \frac{1}{2}$, we have
\begin{equation}
   \label{eq:4.1.2}
\hbox{diam}(\Omega) \leq M_1 r_0,
\end{equation}
\begin{equation}
   \label{eq:4.1.3}
    \hbox{dist}(D, \partial \Omega) \geq 2M_0r_0,
\end{equation}
\begin{equation}
   \label{eq:4.1.4}
\partial\Omega \hbox{ is of } class\  C^{1,\alpha}
\ with\ constants\ r_0, M_0,
\end{equation}
\begin{equation}
   \label{eq:4.1.5}
\partial D \hbox{ is of } class\  C^{6,\alpha}
\ with\ constants\ r_0, M_0,
\end{equation}
with $\alpha$ such that $0<\alpha\leq 1$.

Let us denote by $\Sigma$ the open portion of $\partial \Omega$
where measurements are taken. We assume that there exists
$P_0\in\Sigma$ such that
\begin{equation}
   \label{eq:4.1.6}
   \partial\Omega\cap
R_{r_0,2M_0r_0}(P_0)\subset\Sigma,
\end{equation}
and
\begin{equation}
   \label{eq:4.1.7}
\Sigma \hbox{ is of } class\  C^{2,\alpha} \ with\ constants\
r_0, M_0.
\end{equation}

Let us notice that, without loss of generality, we have chosen
$M_0\geq\frac{1}{2}$ to ensure that $B_{r_0}(P) \subset
R_{r_0,2M_0r_0}(P)$ for every $P \in
\partial \Omega$.

\medskip
\noindent {\it ii) A priori information on the Neumann boundary data.}

We assume that
\begin{equation}
   \label{eq:4.2.1}
\widehat{N}\in H^{-\frac{1}{2}}(\partial \Omega,\R^2),\quad
\widehat{N}\not\equiv 0,
\end{equation}
\begin{equation}
  \label{eq:4.2.2}
  \int_{\partial\Omega}\widehat{N}\cdot r =0, \quad \hbox{for every } r\in {\mathcal R}_2,
\end{equation}
\begin{equation}
   \label{eq:4.2.2bis}
\hbox{supp}(\widehat{N})\subset\subset\Sigma,
\end{equation}
and that, for a given constant $F>0$,
\begin{equation}
\label{eq:4.2.3}
   \frac{\|\widehat{N}\|_{H^{-\frac{1}{2}}(\partial
   \Omega ,\R^2)}}
    {\|\widehat{N}\|_{H^{-1}(\partial \Omega,\R^2)}}\leq
   F,
\end{equation}

{\it iii) A  priori information on the elasticity tensor.}

The constitutive equation \eqref{eq:I4} can be written as
\begin{equation}
\label{eq:4.2.4}
   N_{\alpha\beta}(x)=C_{\alpha\beta\gamma\delta}(x)\epsilon_{\gamma\delta},
\end{equation}
where the elasticity tensor $\C=(C_{\alpha\beta\gamma\delta})$ is defined as
\begin{equation}
  \label{eq:4.2.5}
   \C (x) A = \frac{Eh}{1-\nu^2(x)} ((1-\nu(x))\widehat{A}+ \nu(\hbox{tr}(A))I_2),
\end{equation}
for every $2 \times 2$ matrix $A$, where the Young's modulus $E$ and the Poisson's coefficient $\nu$ are given in terms of the Lam\'{e} moduli
as follows
\begin{equation}
  \label{eq:E_nu}
  E(x)=\frac{\mu(x)(2\mu(x)+3\lambda(x))}{\mu(x)+\lambda(x)},\qquad\nu(x)=\frac{\lambda(x)}{2(\mu(x)+\lambda(x))}.
\end{equation}

On the Lam\'e coefficients $\mu=\mu(x)$, $\lambda=\lambda(x)$, $\mu:\overline{\Omega}\rightarrow\R$, $\lambda:\overline{\Omega}\rightarrow\R$, we assume
\begin{equation}
  \label{eq:4.3.1}
  \mu(x)\geq \alpha_0,\qquad 2\mu(x)+3\lambda(x)\geq\gamma_0, \qquad \hbox{in } \overline{\Omega},
\end{equation}
for positive constants $\alpha_0$ and $\gamma_0$.

The above assumptions ensure that $\C$ satisfies the minor and
major symmetries
\begin{equation}
  \label{eq:minor-major-symm}
C_{\alpha\beta\gamma\delta}=C_{\beta\alpha\gamma\delta}=C_{\alpha\beta\delta\gamma},
\ \ C_{\alpha\beta\gamma\delta}=C_{\gamma\delta\alpha\beta}, \quad
\hbox{for every } \alpha, \beta, \gamma, \delta=1,2, \ \hbox{in
}\overline{\Omega},
\end{equation}
and that it is strongly convex in
$\overline{\Omega}$, precisely
\begin{equation}
  \label{eq:strong_convexity}
  \mathbb{C} A\cdot A\geq h\xi_0|A|^2, \qquad \hbox{in }\Omega,
\end{equation}
for every $2\times 2$ symmetric matrix $A$, where
$\xi_0=\min\{2\alpha_0,\gamma_0\}$ (see \cite[Lemma 3.5]{INDIANA} for details).
Moreover, $E(x)>0$ and $-1<\nu(x)<\frac{1}{2}$ in
$\overline{\Omega}$.

We further assume that
\begin{equation}
  \label{eq:4.3.5}
  \|\lambda\|_{C^4(\overline{\Omega})}, \quad \|\mu\|_{C^4(\overline{\Omega})}\leq \Lambda_0,
\end{equation}
for some positive constant $\Lambda_0$.

We note that the equilibrium problem \eqref{eq:I1}--\eqref{eq:I5} can be written in compact form as
\begin{center}
\( {\displaystyle \left\{
\begin{array}{lr}
     \divrg(\C \nabla a)=0,
      & \mathrm{in}\ \Omega \setminus \overline {D},
        \vspace{0.25em}\\
      (\C\nabla a)n=\widehat{N}, & \mathrm{on}\ \partial \Omega,
          \vspace{0.25em}\\
      (\C\nabla a)n=0, & \mathrm{on}\ \partial D.
        \vspace{0.25em}\\
\end{array}
\right. } \) \vskip -5.9em
\begin{eqnarray}
& & \label{eq:4.3.2}\\
& & \label{eq:4.3.3}\\
& & \label{eq:4.3.4}
\end{eqnarray}
\end{center}
The weak formulation of \eqref{eq:4.3.2}--\eqref{eq:4.3.4}
consists in finding $a=a(x)\in H^1(\Omega \setminus \overline
{D})$ satisfying
\begin{equation}
  \label{eq:4.3.6}
  \int_{\Omega \setminus \overline {D}}\C\nabla a\cdot \nabla v =\int_{\partial\Omega}\widehat{N}\cdot v, \quad \hbox{for every }v\in H^1(\Omega \setminus \overline {D}).
\end{equation}
Under our assumptions, there exists a unique solution to
\eqref{eq:4.3.6} up to addition of a rigid displacement. In order
to select a single solution, we shall assume the normalization
conditions
\begin{equation}
  \label{eq:4.4.1}
  \int_{\Omega \setminus \overline {D}}a =0, \quad \int_{\Omega \setminus \overline {D}}(\nabla a - \nabla^Ta) =0,
\end{equation}
which imply the following stability estimate for the direct
problem \eqref{eq:4.3.2}--\eqref{eq:4.3.4}
\begin{equation}
  \label{eq:4.4.2}
  \|a\|_{H^1(\Omega \setminus \overline {D})}\leq Cr_0\|\widehat{N}\|_{H^{-\frac{1}{2}}(\partial\Omega,\R^2)},
\end{equation}
where $C>0$ is a constant only depending on $h$, $\alpha_0$,
$\gamma_0$, $M_0$ and $M_1$.

In what follows, we shall refer to the set of constants $h$,
$\alpha_0$, $\gamma_0$, $\Lambda_0$, $\alpha$, $M_0$, $M_1$ and
$F$ as the a priori data.

\begin{theo}[Stability result]
  \label{theo:Main}
Let $\Omega$ be a domain satisfying \eqref{eq:4.1.2},
\eqref{eq:4.1.4} and let $\Sigma$ be an open portion of
$\partial\Omega$ satisfying \eqref{eq:4.1.6}--\eqref{eq:4.1.7}.
Let the elasticity tensor $\C=\C(x)\in {\mathcal L}(\M^2,\M^2)$
given by \eqref{eq:4.2.5}, with Lam\'e moduli
$\lambda=\lambda(x)$, $\mu=\mu(x)$ satisfying \eqref{eq:4.3.1} and
\eqref{eq:4.3.5}. Let $\widehat{N}\in
H^{-\frac{1}{2}}(\partial\Omega,\R^2)$, $\widehat{N}\not\equiv 0$,
satisfying \eqref{eq:4.2.2}--\eqref{eq:4.2.3}. Let $D_i$, $i=1,2$,
be two open subsets of $\Omega$ satisfying \eqref{eq:4.1.1},
\eqref{eq:4.1.3}, \eqref{eq:4.1.5}, and let $a^{(i)}\in H^1(\Omega
\setminus \overline{D_i}, \R^2)$ be the solution to
\eqref{eq:4.3.2}--\eqref{eq:4.3.4}, satisfying \eqref{eq:4.4.1},
when $D=D_i$, $i=1,2$. If, given $\varepsilon>0$, we have
\begin{equation}
    \label{eq:4.5.1}
\min_{r \in {\mathcal R}_2} \|a^{(1)} - a^{(2)} -r \|_{L^2(\Sigma, \R^2)}\leq
r_0\varepsilon,
\end{equation}
then we have
\begin{equation}
    \label{eq:4.5.2}
d_{\cal H}( \overline{D_1},\overline{D_2} ) \leq
Cr_0\left|\log\left(\frac{\varepsilon}{\|\widehat{N}\|_{H^{-\frac{1}{2}}(\partial\Omega,\R^2)}}\right)\right|^{-\eta},
\end{equation}
where $C$, $\eta$, $C>0$, $\eta>0$, only
depend on the a priori data.
\end{theo}

\begin{rem}
\label{rem:multiple}
Let us notice that, as it will be clear {}from the proof, the
above stability result holds true also when the domain $\Omega$
contains a finite number of connected cavities $D^{(j)}$, $j=1,
\dots, J$, such that $\partial D^{(j)} \in C^{6,\alpha}$ with
constants $r_0$, $M_0$, and $\hbox{dist}( \partial D^{(j)},
\partial D^{(k)}) \geq r_0$, for $j \neq k$.
\end{rem}

\section{Airy's transformation}
   \label{sec:Airy}

It is known that the boundary value problem in plane linear
elasticity can be formulated in terms of an equivalent
Kirchhoff-Love \textit{plate-like} problem involving a
scalar-valued function called \textit{Airy's function}. Although
this argument is well established, see, for instance,
\cite{Gurtin} and \cite{Fichera}, for reader convenience in what
follows we recall the essential points of the analysis.

For the sake of completeness, we consider a mixed boundary value
problem, in order to describe the transformation of both Dirichlet
and Neumann boundary conditions. Let $a=a_1e_1+a_2e_2$, $a \in
H^1(\mathcal U, \R^2)$, be the solution to the GPS problem
\begin{center}
\( {\displaystyle \left\{
\begin{array}{lr}
     N_{\alpha\beta,\beta}=0,
      & \mathrm{in}\ \mathcal U,
        \vspace{0.25em}\\
      N_{\alpha\beta}n_\beta=\widehat{N}_\alpha, & \mathrm{on}\ \partial_t \mathcal U,
          \vspace{0.25em}\\
      a_\alpha = \widehat{a}_\alpha, & \mathrm{on}\ \partial_u \mathcal U,
        \vspace{0.25em}\\
      N_{\alpha\beta}=\frac{Eh}{1-\nu^2}\left((1-\nu)\epsilon_{\alpha\beta}+\nu(\epsilon_{\gamma\gamma})\delta_{\alpha\beta}\right), &\mathrm{in}\ \mathcal U,
          \vspace{0.25em}\\
        \epsilon_{\alpha\beta}=\frac{1}{2}(a_{\alpha,\beta}+a_{\beta,\alpha}), &\mathrm{in}\
        \mathcal U,
          \vspace{0.25em}\\
\end{array}
\right. } \) \vskip -8.9em
\begin{eqnarray}
& & \label{eq:AIRY1-1}\\
& & \label{eq:AIRY1-2}\\
& & \label{eq:AIRY1-3}\\
& & \label{eq:AIRY1-4}\\
& & \label{eq:AIRY1-5}
\end{eqnarray}
\end{center}
where $\widehat{N} \in H^{-1/2}(\partial_t \mathcal U, \R^2)$ and
$\widehat{a} \in H^{1/2}(\partial_u \mathcal U, \R^2)$ are given
Neumann and Dirichlet data, respectively. Here, $\partial_u
\mathcal U$, $\partial_t \mathcal U$ are two disjoint connected open
subsets of $\partial \mathcal U$, with $\partial \mathcal U =
\overline{\partial_u \mathcal U \cup \partial_t \mathcal U}$.

The equilibrium equations \eqref{eq:AIRY1-1}, and the simply
connectness of $\mathcal U$, ensure the existence of a single-valued
function $\varphi=\varphi(x_1,x_2)$, $\varphi \in H^2(\mathcal U)$,
such that
\begin{equation}
   \label{eq:AIRY2-1}
   N_{\alpha \beta}= e_{\alpha \gamma} e_{\beta \delta}
   \varphi_{,\gamma \delta},
\end{equation}
where the matrix $e_{\alpha \gamma}$ is defined as follows:
$e_{11}=e_{22}=0$, $e_{12}=1$, $e_{21}=-1$; see \cite{Airy}. We
recall that, by construction, the function $\varphi$ and its first
partial derivatives $\varphi_{,1}$, $\varphi_{,2}$ are uniquely
determined up to an additive arbitrary constant.

It is convenient to introduce the \textit{strain functions}
$K_{\alpha \beta}$, $\alpha,\beta=1,2$, associated to the
infinitesimal strain $\epsilon_{\alpha \beta}$:
\begin{equation}
   \label{eq:AIRY2-2}
   K_{\alpha \beta}= e_{\delta\alpha } e_{\gamma\beta}
   \epsilon_{\delta \gamma }, \quad \alpha,\beta=1,2.
\end{equation}
By inverting the constitutive equation \eqref{eq:AIRY1-4}, we get
\begin{equation}
   \label{eq:AIRY2-2bis}
   \epsilon_{\alpha\beta}=\frac{1+\nu}{Eh} N_{\alpha \beta} -
   \frac{\nu}{Eh}(N_{\gamma \gamma}) \delta_{\alpha\beta},
\end{equation}
and using \eqref{eq:AIRY2-1} we obtain
\begin{equation}
   \label{eq:AIRY2-2ter}
   \epsilon_{\alpha\beta}=\frac{1+\nu}{Eh} e_{\alpha \gamma} e_{\beta \delta}
   \varphi_{,\gamma \delta} -
   \frac{\nu}{Eh}(\varphi_{,\gamma \gamma}) \delta_{\alpha\beta}.
\end{equation}
Inserting this expression of $\epsilon_{\alpha\beta}$ into
\eqref{eq:AIRY2-2}, we have
\begin{equation}
   \label{eq:AIRY2-3}
   K_{\alpha \beta}= L_{\alpha \beta \gamma \delta} \varphi_{,\gamma
   \delta},
\end{equation}
where the Cartesian components $L_{\alpha \beta \gamma \delta}$ of
the fourth order tensor $\mathbb{L}$ are
\begin{equation}
   \label{eq:AIRY2-4}
   L_{\alpha \beta \gamma \delta}=
   \frac{1+\nu}{Eh} \delta_{\alpha \gamma} \delta_{\beta \delta}
   - \frac{\nu}{Eh} \delta_{\alpha \beta} \delta_{\gamma \delta}.
\end{equation}
The strain $\epsilon_{\alpha \beta}$ obviously satisfies the
well-known two-dimensional \textit{Saint-Venant compatibility
equation}
\begin{equation}
   \label{eq:AIRY2-5}
   \epsilon_{11,22}+\epsilon_{22,11}-2\epsilon_{12,12}=0, \quad \hbox{in }
   {\mathcal U}.
\end{equation}
Inverting \eqref{eq:AIRY2-2}, we have
\begin{equation}
   \label{eq:AIRY3-1}
   \epsilon_{\alpha \beta}=e_{\alpha \gamma} e_{\beta \delta} K_{\gamma \delta} ,
\end{equation}
and the equation \eqref{eq:AIRY2-5}, written in terms of
$K_{\gamma \delta}$, becomes
\begin{equation}
   \label{eq:AIRY3-2}
   {\divrg} ( {\divrg} ( \mathbb{L} \nabla^2 \varphi))=0, \quad \hbox{in }
   {\mathcal U},
\end{equation}
or, more explicitly,
\begin{multline}
   \label{eq:AIRY3-3}
\Delta^2\varphi+2Eh \nabla \left(\frac{1}{Eh} \right ) \cdot
\nabla(\Delta\varphi)
-Eh\Delta\left(\frac{\nu}{Eh}\right)\Delta\varphi+\\
+Eh \nabla^2\left(\frac{1+\nu}{Eh}\right) \cdot \nabla^2 \varphi
=0,\quad \hbox{in } \mathcal U.
\end{multline}
The above partial differential equation expresses the form assumed
by the field equation \eqref{eq:AIRY1-1} in terms of the Airy's
function $\varphi$.

We now consider the transformation of the Neumann boundary
condition \eqref{eq:AIRY1-2} on $\partial_t \mathcal U$. By
\eqref{eq:AIRY2-1}, the condition on $\partial_t \mathcal U$ can be
written as
\begin{equation}
   \label{eq:AIRY3-4}
   e_{\alpha \gamma} e_{\beta \delta} \varphi_{,\gamma \delta}
   n_\beta =\widehat{N}_{\alpha},
\end{equation}
that is, recalling that $\tau_\delta = e_{\beta \delta} n_\beta$
on $\partial \mathcal U$,
\begin{equation}
   \label{eq:AIRY3-5}
   (\varphi_{,1})_{,s}= - \widehat{N}_2, \quad (\varphi_{,2})_{,s}=
   \widehat{N}_1, \quad \hbox{on } \partial_t \mathcal U,
\end{equation}
where $s$ is an arc length parametrization on $\partial \mathcal
U$. By integrating the above equations with respect to $s$, {}from
$P_0 \in
\partial_t \mathcal U$ to $P \in \partial_t \mathcal U$, with $s(P_0)=0$
and $s(P)=s$, the gradient of $\varphi$ on $\partial_t \mathcal U$ can
be determined up to an additive constant vector $c=c_1 e_1 + c_2
e_2$, namely
\begin{equation}
   \label{eq:AIRY4-1}
   \nabla \varphi (s) = c + \widehat{g}(s), \quad \hbox{on }
   \partial_t \mathcal U,
\end{equation}
where $\widehat{g}(s) =\widehat{g}_1(s)e_1 + \widehat{g}_2(s)e_2$,
$\widehat{g}_1(s) = - \int_0^s \widehat{N}_2(\xi)d\xi$,
$\widehat{g}_2(s) = \int_0^s \widehat{N}_1(\xi)d\xi$. It follows
that the normal derivative of $\varphi$ on $\partial_t \mathcal U$ is
prescribed in terms of the Neumann data $\widehat{N}$, that is,
\begin{equation}
   \label{eq:AIRY4-2}
   \varphi_{,n} = (c + \widehat{g}(s))\cdot n, \quad \hbox{on }
   \partial_t \mathcal U,
\end{equation}
whereas, integrating once more \eqref{eq:AIRY4-1} {}from $P_0$ to
$P$, we have
\begin{equation}
   \label{eq:AIRY4-3}
   \varphi(s) = C + \widehat{G}(s), \quad \hbox{on }
   \partial_t \mathcal U,
\end{equation}
where $C=\varphi(0)=$constant, and $\widehat{G}(s) = \int_0^s (c +
\widehat{g}(\xi)) \cdot \tau(\xi) d\xi$. We notice that it is
always possible to select the two arbitrary constants occurring in
the construction of $\nabla \varphi$ such that $c_1=c_2=0$ (see,
for example, \cite{Sokolnikoff} for details). In particular, if
the Neumann data $\widehat{N}$ vanishes on $\partial_t \mathcal
U$, then we can also choose the third constant $C=0$, so that
$\varphi(s)=0$ on $\partial_t \mathcal U$. In this case, the
homogeneous Neumann boundary conditions for the GPS problem are
transformed into the homogeneous Dirichlet boundary conditions for
the Airy's function:
\begin{equation}
   \label{eq:AIRY4-4}
   \varphi=0, \quad  \varphi_{,n}=0, \quad \hbox{on }
   \partial_t \mathcal U.
\end{equation}
The determination of the boundary conditions satisfied by
$\varphi$ on $\partial_u \mathcal U$ is less obvious, since the
corresponding boundary conditions in the original two-dimensional
elasticity problem are not explicitly expressed in terms of the
Airy's function or its derivatives. In dealing with this boundary
condition, we need to assume $C^{1,1}$-regularity for $\partial
\mathcal U$. We adopt a variational-like approach. Without loss of
generality, we can assume $\partial_u \mathcal U =
\partial \mathcal U$.

Let $\widetilde{\varphi}$, $\widetilde{\varphi}: \overline{\mathcal U}
\rightarrow \R$, be a $C^\infty$-test function, and define the
\textit{associated Airy stress field}
\begin{equation}
   \label{eq:AIRY5-1}
   \widetilde{N}_{\alpha \beta}= e_{\alpha \gamma} e_{\beta \delta}
   \widetilde{\varphi}_{, \gamma \delta}, \quad \hbox{in }
   \overline{\mathcal U},
\end{equation}
which obviously satisfies the equilibrium equations
\begin{equation}
   \label{eq:AIRY5-2}
   \widetilde{N}_{\alpha \beta, \beta}= 0, \quad \hbox{in }
   {\mathcal U}.
\end{equation}
Multiplying \eqref{eq:AIRY5-2} by the displacement field $a=a_1
e_1 + a_2 e_2$ solution to \eqref{eq:AIRY1-1}--\eqref{eq:AIRY1-5},
and integrating by parts, we obtain
\begin{equation}
   \label{eq:AIRY5-3}
   \int_{\mathcal U} \widetilde{\varphi}_{, \gamma \delta}  K_{\gamma \delta}=
   \int_{\partial {\mathcal U}} \widetilde{N}_{\alpha \beta} n_\beta \widehat{a}_\alpha.
\end{equation}
We first work on the integral of the left hand side of
\eqref{eq:AIRY5-3}. After two integrations by parts, we obtain
\begin{equation}
   \label{eq:AIRY6-1}
   \int_{\mathcal U} \widetilde{\varphi}_{, \gamma \delta}  K_{\gamma \delta}=
   \int_{\mathcal U} K_{\gamma \delta, \gamma \delta} \widetilde{\varphi}
   +
   \int_{\partial {\mathcal U}} \widetilde{\varphi}_{,\gamma} K_{\gamma
   \delta}n_\delta
   -
   \int_{\partial {\mathcal U}}
   \widetilde{\varphi} K_{\gamma
   \delta, \delta} n_\gamma.
\end{equation}
We elaborate the second integral $I$ on the right hand side of the above equation in terms of the local
coordinates. Recalling that $\tau_\alpha = e_{\beta
\alpha}n_\beta$ on $\partial {\mathcal U}$ and $
\widetilde{\varphi}_{,\alpha} = n_\alpha \widetilde{\varphi}_{,n}
+\tau_\alpha \widetilde{\varphi}_{,s}$ on $\partial {\mathcal U}$,
$\alpha, \beta =1,2$, we have
\begin{equation}
   \label{eq:AIRY6-2}
   I = \int_{\partial {\mathcal U}}
   ( \widetilde{\varphi}_{,n} K_{nn} + \widetilde{\varphi}_{,s}
   K_{\tau n}
   ),
\end{equation}
where, to simplify the notation, we have introduced on $\partial
{\mathcal U}$ the two functions
\begin{equation}
   \label{eq:AIRY6-3}
   K_{nn}= K_{\gamma \delta}n_\delta n_\gamma, \quad K_{n \tau} =
   K_{\gamma \delta}n_\delta \tau_\gamma (=K_{\tau n}).
\end{equation}
Since $\partial {\mathcal U}$ is of $C^{1,1}$-class, integrating by
parts the second term in \eqref{eq:AIRY6-2} gives
\begin{equation}
   \label{eq:AIRY6-4}
   I = \int_{\partial {\mathcal U}}
   ( \widetilde{\varphi}_{,n} K_{nn} - \widetilde{\varphi}
   K_{\tau n, s}
   ).
\end{equation}
Therefore, the left hand side of \eqref{eq:AIRY5-3} takes the form
\begin{equation}
   \label{eq:AIRY6-5}
   \int_{\mathcal U} \widetilde{\varphi}_{, \gamma \delta}  K_{\gamma \delta}=
   \int_{\mathcal U} K_{\gamma \delta, \gamma \delta} \widetilde{\varphi}
   +
   \int_{\partial {\mathcal U}} ( K_{nn} \widetilde{\varphi}_{,n}
   -
    (K_{\gamma
   \delta, \delta} n_\gamma
    +
    K_{\tau n, s} ) \widetilde{\varphi} ).
\end{equation}
We next elaborate the integral appearing on the right hand side of
\eqref{eq:AIRY5-3}. Let us introduce the \textit{boundary
displacement functions} associated to the Dirichlet data
$\widehat{a}$:
\begin{equation}
   \label{eq:AIRY7-1}
   \widehat{U}_\gamma = e_{\alpha \gamma}\widehat{a}_\alpha, \quad
   \hbox{on } \partial {\mathcal U}.
\end{equation}
Passing to local coordinates, after an integration by parts, we
have
\begin{equation}
   \label{eq:AIRY7-2}
   \int_{\partial {\mathcal U}} \widetilde{N}_{\alpha \beta} n_\beta \widehat{a}_\alpha
   =
   \int_{\partial {\mathcal U}} \widetilde{\varphi}_{,\gamma \delta}
   \tau_\delta \widehat{U}_\gamma=
   \int_{\partial {\mathcal U}} (\widetilde{\varphi}_{,\gamma})_{,s} \widehat{U}_\gamma
   =
   - \int_{\partial {\mathcal U}} \widetilde{\varphi}_{,\gamma}
    \widehat{U}_{\gamma,s}.
\end{equation}
Expressing again $\nabla \widetilde{\varphi}$ in terms of local
coordinates, and integrating by parts, by the regularity of
$\partial {\mathcal U}$ we obtain
\begin{equation}
   \label{eq:AIRY7-3}
   \int_{\partial {\mathcal U}} \widetilde{N}_{\alpha \beta} n_\beta \widehat{a}_\alpha
   =
   \int_{\partial {\mathcal U}} (- \widetilde{\varphi}_{n} \widehat{U}_{\gamma,s}n_\gamma
   +
   \widetilde{\varphi} (\tau_\gamma \widehat{U}_{\gamma,s})_{,s}).
\end{equation}
Finally, by rewriting \eqref{eq:AIRY5-3} using \eqref{eq:AIRY6-5}
and \eqref{eq:AIRY7-3}, the strain functions $K_{\gamma \delta}$
satisfy the condition
\begin{equation}
   \label{eq:AIRY7-4}
   \int_{\mathcal U} K_{\gamma \delta, \gamma \delta} \widetilde{\varphi}
   +
   \int_{\partial {\mathcal U}} ( K_{nn} + \widehat{U}_{\gamma,s}n_\gamma ) \widetilde{\varphi}_{,n}
   -
   \int_{\partial {\mathcal U}}  (K_{\gamma
   \delta, \delta} n_\gamma
    +
    K_{\tau n, s}  + (\tau_\gamma \widehat{U}_{\gamma,s})_{,s} )
    \widetilde{\varphi} =0,
\end{equation}
for every $\widetilde{\varphi} \in C^\infty (\overline{{\mathcal U}})$.
By the arbitrariness of the test function $\widetilde{\varphi}$,
and of the traces of $\widetilde{\varphi}$ and
$\widetilde{\varphi}_{,n}$ on $\partial {\mathcal U}$, we determine the
conditions satisfied by $K_{\gamma \delta}$, namely, the field
equation
\begin{equation}
   \label{eq:AIRY8-1}
   K_{\gamma \delta, \gamma \delta} =0, \quad \hbox{in } {\mathcal U},
\end{equation}
which coincides with \eqref{eq:AIRY3-2}, and the two boundary
conditions
\begin{equation}
   \label{eq:AIRY8-2}
   K_{nn} = - \widehat{U}_{\gamma,s}n_\gamma, \quad \hbox{on } \partial {\mathcal U},
\end{equation}
\begin{equation}
   \label{eq:AIRY8-3}
   K_{\gamma\delta, \delta} n_\gamma
    +
    K_{\tau n, s}  = - (\tau_\gamma \widehat{U}_{\gamma,s})_{,s}, \quad \hbox{on } \partial
    {\mathcal U}.
\end{equation}
The above equations \eqref{eq:AIRY8-1} and \eqref{eq:AIRY8-2},
\eqref{eq:AIRY8-3} are known as \textit{compatibility field
equation} and \textit{compatibility boundary conditions} for the
strain functions $K_{\gamma \delta}$, respectively. In conclusion,
under the assumption $\widehat{N}=0$ on $\partial_t {\mathcal U}$,
the two-dimensional elasticity problem
\eqref{eq:AIRY1-1}--\eqref{eq:AIRY1-5} can be formulated in terms
of the Airy's function as follows:
\begin{center}
\( {\displaystyle \left\{
\begin{array}{lr}
     K_{\gamma \delta, \gamma \delta} =0,
      & \mathrm{in}\ {\mathcal U},
        \vspace{0.25em}\\
      \varphi=0, & \mathrm{on}\ \partial_t {\mathcal U},
          \vspace{0.25em}\\
      \frac{\partial \varphi}{\partial n}=0, & \mathrm{on}\ \partial_t {\mathcal U},
        \vspace{0.25em}\\
      K_{\alpha \beta} n_{\alpha} n_{\beta}= - \widehat{U}_{\gamma,s}n_\gamma, & \mathrm{on}\ \partial_u {\mathcal U},
        \vspace{0.25em}\\
      K_{\alpha\beta, \beta} n_\alpha + (K_{\alpha \beta} n_\beta \tau_\alpha)_{,s}  =
       - (\tau_\gamma \widehat{U}_{\gamma,s})_{,s}, & \mathrm{on}\ \partial_u {\mathcal U},
        \vspace{0.25em}\\
      K_{\alpha\beta}=\frac{1}{Eh}
      \left(
      (1+\nu)\varphi_{,\alpha\beta}-\nu(\Delta \varphi) \delta_{\alpha\beta}
      \right),
       & \mathrm{in}\ \overline{{\mathcal U}}.
          \vspace{0.25em}\\
\end{array}
\right. } \) \vskip -10.3em
\begin{eqnarray}
& & \label{eq:AIRY8-4}\\
& & \label{eq:AIRY8-5}\\
& & \label{eq:AIRY8-6}\\
& & \label{eq:AIRY8-7}\\
& & \label{eq:AIRY8-8}\\
& & \label{eq:AIRY8-9}
\end{eqnarray}
\end{center}
There is an important analogy connected with the above boundary
value problem. Equations \eqref{eq:AIRY8-4}--\eqref{eq:AIRY8-9}
describe the conditions satisfied by the transversal displacement
$\varphi=\varphi(x_1,x_2)$ of the middle surface ${\mathcal U}$ of
a Kirchhoff-Love thin elastic plate made by isotropic material.
The plate is clamped on $\partial_t {\mathcal U}$, and subject to
a couple field $\widehat{M}=\widehat{M}_\tau n + \widehat{M}_n
\tau$ assigned on $\partial_u {\mathcal U}$, with $\widehat{M}_n=
- \widehat{U}_{\gamma,s}n_\gamma$ and $\widehat{M}_\tau=
\tau_\gamma \widehat{U}_{\gamma,s}$, see, for example,
\cite{INDIANA}. Within this analogy, the strain functions
$K_{\alpha\beta}=K_{\alpha\beta}(x_1,x_2)$ play the role of the
\textit{bending moments} (for $\alpha=\beta$) and the
\textit{twisting moments} (for $\alpha \neq \beta$) of the plate
at $(x_1, x_2) \in \overline{\Omega}$ (per unit length), and the
\textit{bending stiffness} of the plate is equal to $(Eh)^{-1}$.

Let us observe that the geometry of the inverse problem here
considered, that is $\mathcal U = \Omega \setminus \overline{D}$
does not ensure the existence of a globally defined Airy's
function, since the hypotheses of simple connectedness is missing.
For this reason, in the following Section \ref{sec:proof} we shall
make use of local Airy's functions, defined either in interior
discs (see the proof of Proposition \ref{prop:FVRI}) or in
neighbourhoods of the boundary of the cavity (see the proof of
Proposition \ref{prop:FVRB}).

\begin{prop}
  \label{prop:15_3}
Under the above notation and assumptions, we have
\begin{equation}
   \label{eq:15.3}
\frac{(1-|\nu|)^2}{E^2h^2}|\nabla^2
\varphi|^2\leq|\widehat{\nabla} a|^2\leq
\frac{(1+|\nu|)^2}{E^2h^2}|\nabla^2 \varphi|^2
\end{equation}

\end{prop}

\begin{proof}

By \eqref{eq:AIRY2-1}, we have $N_{11} = \varphi_{,22}$, $N_{22} =
\varphi_{,11}$, $N_{12} = N_{21} = -\varphi_{,12}$, so that

\begin{equation}
   \label{eq:Edi15.1}
|\nabla^2
\varphi|^2=\sum_{\alpha,\beta=1}^2N_{\alpha\beta}^2.
\end{equation}

By \eqref{eq:AIRY2-2bis}, we have $\epsilon_{11}=\frac{1}{Eh}
N_{11}  -
   \frac{\nu}{Eh}N_{22}$, $\epsilon_{22}=\frac{1}{Eh} N_{22}  -
   \frac{\nu}{Eh}N_{11}$, $\epsilon_{12}=\epsilon_{21}=\frac{1+\nu}{Eh} N_{12}$, so that

\begin{equation}
   \label{eq:Edi15.2}
|\widehat{\nabla}
a|^2=\sum_{\alpha,\beta=1}^2\epsilon_{\alpha\beta}^2=\frac{1}{(Eh)^2}\{(1+\nu^2)(N_{11}^2 +N_{22}^2)+2(1+\nu)^2N_{12}^2-4\nu N_{11}N_{22}\}.
\end{equation}
Let us estimate the term $-4\nu N_{11}N_{22}$ by using the elementary inequalities
\begin{equation}
   \label{eq:Edi_elementary}
\pm 2N_{11}N_{22}\leq N_{11}^2+N_{22}^2.
\end{equation}

\emph{I)} Estimate {}from below.

\emph{i)} $0<\nu<\frac{1}{2}$

If $N_{11}N_{22}<0$, then $-4\nu N_{11}N_{22}>0$, whereas if
$N_{11}N_{22}\geq 0$, then, by \eqref{eq:Edi_elementary}, $-4\nu
N_{11}N_{22}\geq -2\nu(N_{11}^2+N_{22}^2)$. Since
$-2\nu(N_{11}^2+N_{22}^2)\leq 0$, we have, independently of the
sign of $N_{11}N_{22}$,

\begin{equation}
   \label{eq:Edi_basso_nu_positivo}
    -4\nu N_{11}N_{22}\geq -2\nu(N_{11}^2+N_{22}^2).
\end{equation}

\emph{ii)} $\nu=0$

In  this case,

\begin{equation}
   \label{eq:Edi_basso_nu_nullo}
    -4\nu N_{11}N_{22}=0.
\end{equation}

\emph{iii)} $-1<\nu<0$ ($\Leftrightarrow 0<-\nu<1$)

If $N_{11}N_{22}\geq 0$, then $-4\nu N_{11}N_{22}\geq 0$, whereas
if $N_{11}N_{22}< 0$, then, by \eqref{eq:Edi_elementary}, $-4\nu
N_{11}N_{22}\geq 2\nu(N_{11}^2+N_{22}^2)$. Since
$2\nu(N_{11}^2+N_{22}^2)\leq 0$, we have, independently of the
sign of $N_{11}N_{22}$,

\begin{equation}
   \label{eq:Edi_basso_nu_negativo}
    -4\nu N_{11}N_{22}\geq 2\nu(N_{11}^2+N_{22}^2).
\end{equation}

Therefore, collecting together the three cases, we have
\begin{equation}
   \label{eq:Edi15.4}
-4\nu N_{11}N_{22}\geq -2|\nu|(N_{11}^2+N_{22}^2).
\end{equation}

{}From \eqref{eq:Edi15.2} and \eqref{eq:Edi15.4}, we have

\begin{multline}
   \label{eq:Edi16.1}
|\widehat{\nabla}
a|^2\geq
\frac{1}{(Eh)^2}\{(1+|\nu|^2-2|\nu|)(N_{11}^2 +N_{22}^2)+2(1+\nu)^2N_{12}^2\}
\geq \\
\geq\frac{(1-|\nu|)^2}{E^2h^2}\sum_{\alpha,\beta =1}^2
N_{\alpha\beta}^2= \frac{(1-|\nu|)^2}{E^2h^2}|\nabla^2 \varphi|^2.
\end{multline}

\emph{II)} Estimate {}from above.

By distinguishing the three cases as above, we get similarly
\begin{equation}
   \label{eq:Edi16.2}
-4\nu N_{11}N_{22}\leq 2|\nu|(N_{11}^2+N_{22}^2).
\end{equation}

{}From \eqref{eq:Edi15.2} and \eqref{eq:Edi16.2}, we get the right hand side
of \eqref{eq:15.3}.

\end{proof}

\section{Proof of the main result}
   \label{sec:proof}

\begin{prop}[Lipschitz Propagation of Smallness]
  \label{prop:LPS}
Let $\Omega$ be a domain satisfying \eqref{eq:4.1.2}, \eqref{eq:4.1.4}. Let $D$ be an open
 subset of $\Omega$ satisfying
\eqref{eq:4.1.1}, \eqref{eq:4.1.3}, \eqref{eq:4.1.5}. Let $a\in
H^1(\Omega \setminus \overline{D}, \R^2)$ be the solution to
\eqref{eq:4.3.2}--\eqref{eq:4.3.4}, satisfying \eqref{eq:4.4.1}.
Let the elasticity tensor $\C=\C(x)\in {\mathcal L}(\M^2,\M^2)$
given by \eqref{eq:4.2.5}, with Lam\'e moduli
$\lambda=\lambda(x)$, $\mu=\mu(x)$ satisfying \eqref{eq:4.3.1} and
\eqref{eq:4.3.5}. Let $\widehat{N}\in
H^{-\frac{1}{2}}(\partial\Omega,\R^2)$, $\widehat{N}\not\equiv 0$,
satisfying \eqref{eq:4.2.2}--\eqref{eq:4.2.3}. Then, there exists
$s>1$, only depending on $\alpha_0$, $\gamma_0$, $\Lambda_0$ and
$M_0$, such that for every $\rho>0$ and every $\bar x\in
(\Omega\setminus \overline{D})_{s\rho}$, we have
\begin{equation}
   \label{eq:LPS}
\int_{B_\rho(\bar x)}|\widehat{\nabla} a|^2\geq
\frac{Cr_0^2}{\exp\left[A\left(\frac{r_0}{\rho}\right)^B\right]}
\|\widehat{N}\|_{H^{-\frac{1}{2}}(\partial\Omega,\R^2)}^2,
\end{equation}
where $A$, $B$, $C>0$ are positive constants only depending on
$\alpha_0$, $\gamma_0$, $\Lambda_0$, $M_0$, $M_1$ and $F$.
\end{prop}

\begin{proof}
The proof follows by merging the Lipschitz Propagation of
Smallness estimate $(3.5)$ contained in \cite[Proposition
3.1]{l:mr2004}, Korn inequalities (see, for instance, \cite{l:f},
\cite{l:amr2008}), trace inequalities (\cite{l:lm}) and
equivalence relations for the $H^{-\frac{1}{2}}$ and
$H^{-1}$-norms of the Neumann data $\widehat{N}$ (see
$(3.9)$--$(3.10)$ in \cite[Remark 3.4]{l:mr2004}).
\end{proof}

\begin{prop}[Finite Vanishing Rate in the Interior]
  \label{prop:FVRI}
Under the hypotheses of Proposition \ref{prop:LPS}, there exist
$\widetilde{c}_0<\frac{1}{2}$ and $C>0$, only depending on
$\alpha_0$, $\gamma_0$ and $\Lambda_0$, such that, for every
$\overline{r}\in(0,r_0)$ and for every $\bar x\in \Omega \setminus
\overline{D}$ such that $B_{\bar r}(\bar x)\subset \Omega
\setminus \overline{D}$, and for every $r_1<\widetilde{c}_0\bar
r$, we have
\begin{equation}
   \label{eq:FVRI}
\int_{B_{r_1}(\bar x)}|\widehat{\nabla} a|^2\geq C
\left(\frac{r_1}{\bar r}\right)^{\tau_0}
\int_{B_{\bar r}(\bar x)}|\widehat{\nabla} a|^2,
\end{equation}
where $\tau_0\geq 1$ only depends on $\alpha_0$,
$\gamma_0$, $\Lambda_0$, $M_0$, $M_1$,$\frac{r_0}{\bar r}$ and $F$.
\end{prop}

\begin{proof}
We can introduce in $B_{\bar r}(\bar x)$ a locally defined Airy's
function $\varphi$ associated to the solution $a$. The proof
follows by adapting the arguments in the proof of the analogous
Proposition 3.5 in \cite{l:mrv2019} which applies to
Kirchhoff-Love plate equation. The main difference consists in
estimating the $L^2$ norms of $\varphi$ and $|\nabla \varphi|$
appearing in $(3.21)$ of \cite{l:mrv2019} in terms of the $L^2$
norm of $|\nabla^2 \varphi|$ and using \eqref{eq:15.3}, the
stability estimate \eqref{eq:4.4.2} and Proposition
\ref{prop:LPS}.

\end{proof}

\begin{prop}[Finite Vanishing Rate at the Boundary]
  \label{prop:FVRB}
Under the hypotheses of Proposition \ref{prop:LPS}, there exist
$\bar c_0<\frac{1}{2}$ and $C>0$, only depending on $\alpha_0$,
$\gamma_0$, $\Lambda_0$, $M_0$, $\alpha$, such that, for every
$\bar x\in \partial D$ and for every $r_1<\bar c_0 r_0$, we have
\begin{equation}
   \label{eq:FVRB}
\int_{B_{r_1}(\bar x)\cap(\Omega\setminus \overline{D})}|\widehat{\nabla} a|^2\geq C
\left(\frac{r_1}{r_0}\right)^{\tau}
\int_{B_{r_0}(\bar x)\cap(\Omega\setminus \overline{D})}|\widehat{\nabla} a|^2,
\end{equation}
where $\tau \geq 1$ only depends on $\alpha_0$,
$\gamma_0$, $\Lambda_0$, $M_0$, $\alpha$, $M_1$ and $F$.
\end{prop}

\begin{proof}
Let us consider the Airy's function $\varphi$ associated to the
solution $a$ and defined in $R_{r_0,2M_0r_0}(\bar x)\cap
\Omega\setminus \overline{D}$, which satisfies the partial
differential equation
\begin{equation}
   \label{eq:AIRY3-div}
   {\divrg} ( {\divrg} ( \mathbb{L} \nabla^2 \varphi))=0, \quad \hbox{in }
   R_{r_0,2M_0r_0}(\bar x)\cap \Omega\setminus
\overline{D},
\end{equation}
or, equivalently,
\begin{multline}
   \label{eq:Ab1.1}
\Delta^2\varphi+2Eh\nabla\left(\frac{1}{Eh}\right ) \cdot
\nabla(\Delta\varphi)
-Eh\Delta\left(\frac{\nu}{Eh}\right)\Delta\varphi+\\
+Eh \nabla^2 \left(\frac{1+\nu}{Eh}\right) \cdot \nabla^2 \varphi
=0,\quad \hbox{in }R_{r_0,2M_0r_0}(\bar x)\cap \Omega\setminus
\overline{D},
\end{multline}
and the homogeneous Dirichlet conditions
\begin{equation}
   \label{eq:Ab1.2}
\varphi=\varphi_{,n}=0,\quad\hbox{ on } \partial D\cap R_{r_0,2M_0r_0}(\bar x).
\end{equation}
Let us notice that, under our assumptions, the fourth order tensor
$\mathbb{L}$ satisfies the strong convexity condition
\begin{equation}
  \label{eq:strong_convexity-L}
  \mathbb{L} A\cdot A\geq  \frac{1}{5h\Lambda_0}  |A|^2, \qquad \hbox{in }\Omega,
\end{equation}
for every $2\times 2$ symmetric matrix $A$. We also notice that
the coefficients of the terms involving second and third-order
derivatives of $\varphi$ in \eqref{eq:Ab1.1} are of class $C^2$ and
$C^3$ in $R_{r_0,2M_0r_0}(\bar x)\cap \Omega\setminus
\overline{D}$, respectively, with corresponding $C^2$ and
$C^3$-norm bounded by a constant only depending on $h$,
$\alpha_0$, $\gamma_0$ and $\Lambda_0$. Therefore, we can apply
the results obtained in \cite{l:arv}. Precisely, by Corollary
$2.3$ in \cite{l:arv}, there exist $c<1$, only depending on $M_0$
and $\alpha$, and $C>1$, only depending on $\alpha_0$, $\gamma_0$,
$\Lambda_0$, $M_0$ and $\alpha$, such that, for every
$r_1<r_2<cr_0$, we have
\begin{equation}
    \label{eq:Ab1.3}
\int_{B_{r_1}(\bar x)\cap(\Omega\setminus \overline{D})}\varphi^2 \geq C\left(\frac{r_1}{r_0}\right)^{\frac{\log B}{\log \frac{cr_0}{r_2}}}
\int_{B_{r_0}(x)\cap(\Omega\setminus \overline{D})}\varphi^2,
\end{equation}
where $B>1$ is given by
\begin{equation}
    \label{eq:Ab1.4}
B= C\left(\frac{r_0}{r_2}\right)^C\frac{\int_{B_{r_0}(\bar x)\cap(\Omega\setminus \overline{D})}\varphi^2}{\int_{B_{r_2}(\bar x)\cap(\Omega\setminus \overline{D})}\varphi^2}.
\end{equation}
Let us choose $r_2 = \overline{c}_0r_0$, with
$\overline{c}_0=\frac{c}{2}$. We need to estimate the quantity
$B$. By applying Poincar\'{e} inequality (see, for instance,
\cite[Example 4.4]{l:amr2008}) and \eqref{eq:15.3}, we have
\begin{equation}
    \label{eq:Ab2.1}
\int_{B_{r_0}(\bar x)\cap(\Omega\setminus \overline{D})}\varphi^2\leq C r_0^4
\int_{B_{r_0}(\bar x)\cap(\Omega\setminus \overline{D})}|\nabla^2\varphi|^2= C r_0^4
\int_{B_{r_0}(\bar x)\cap(\Omega\setminus \overline{D})}|\widehat{\nabla} a|^2
\end{equation}
where $C>0$ only depends on $\alpha_0$, $\gamma_0$, $\Lambda_0$, $M_0$ and $\alpha$.
Moreover, by applying Lemma 4.7 in \cite{l:arv} and \eqref{eq:15.3}, and recalling the choice of $r_2$, we have
\begin{equation}
    \label{eq:Ab2.2}
\int_{B_{r_2}(\bar x)\cap(\Omega\setminus \overline{D})}\varphi^2\geq C r_2^4
\int_{B_{\frac{r_2}{2}}(\bar x)\cap(\Omega\setminus \overline{D})}|\nabla^2\varphi|^2=
C r_0^4\int_{B_{\frac{cr_0}{4}}(\bar x)\cap(\Omega\setminus \overline{D})}|\widehat{\nabla} a|^2.
\end{equation}
By \eqref{eq:Ab2.1}--\eqref{eq:Ab2.2}, using the stability
estimate of the direct problem \eqref{eq:4.4.2} and Proposition
\ref{prop:LPS}, we can estimate $B\leq C$, with $C$ only depending
on $\alpha_0$, $\gamma_0$, $\Lambda_0$, $M_0$, $\alpha$, $M_1$ and
$F$. By using again Poincar\'{e}  inequality, Lemma 4.7 in
\cite{l:arv} and \eqref{eq:15.3}, we obtain the thesis.
\end{proof}

{}From now on, we shall denote by $\mathcal G$ the connected
component of $\Omega \setminus   \overline{(D_1 \cup D_2)}$ such
that $\Sigma \subset \partial{\mathcal G}$.

\begin{prop}[Stability Estimate of Continuation
{from} Cauchy Data]
  \label{prop:4.8.1}
Under the hypotheses of Theorem \ref{theo:Main}, we
have
\begin{equation}
   \label{eq:4.8.1}
\int_{(\Omega\setminus \overline{\mathcal G})\setminus
\overline{D_1}}|\widehat{\nabla} a^{(1)}|^2\leq
r_0^2\|\widehat{N}\|_{H^{-\frac{1}{2}}(\partial \Omega
,\R^2)}^2\omega\left(\frac{\varepsilon}
{\|\widehat{N}\|_{H^{-\frac{1}{2}}(\partial \Omega
,\R^2)}}\right),
\end{equation}
\begin{equation}
   \label{eq:4.8.2}
\int_{(\Omega\setminus \overline{\mathcal G})\setminus
\overline{D_2}}|\widehat{\nabla} a^{(2)}|^2\leq
r_0^2\|\widehat{N}\|_{H^{-\frac{1}{2}}(\partial \Omega
,\R^2)}^2\omega\left(\frac{\varepsilon}
{\|\widehat{N}\|_{H^{-\frac{1}{2}}(\partial \Omega
,\R^2)}}\right),
\end{equation}
where $\omega$ is an increasing continuous function on
$[0,\infty)$ which satisfies
\begin{equation}
   \label{eq:4.8.3}
\omega(t)\leq C(\log|\log t|)^{-\frac{1}{2}},\qquad \hbox{for every }
t<e^{-1},
\end{equation}
with $C>0$ only depending on
$\alpha_0$, $\gamma_0$, $\Lambda_0$, $M_0$, $\alpha$ and $M_1$.
Moreover, there exists $d_0>0$, with $\frac{d_0}{r_0}$ only depending on $M_0$ and $\alpha$,
such that if
$d_{\cal H}( \overline{\Omega\setminus D_1},\overline{\Omega\setminus D_2} )\leq d_0$ then
\eqref{eq:4.8.1}--\eqref{eq:4.8.2} hold with $\omega$
given by
\begin{equation}
   \label{eq:4.8.4}
\omega(t)\leq C|\log t|^{-\sigma},\qquad \hbox{for  every }
t<1,
\end{equation}
where $\sigma>0$ and $C>0$ only depend on $\alpha_0$, $\gamma_0$, $\Lambda_0$,
$M_0$, $\alpha$, $M_1$.
\end{prop}

\begin{proof}
The proof can be easily obtained by adapting the proof of the
analogous estimates contained in Proposition 3.5 and Proposition
3.6 in \cite{l:mr2004}. The only difference consists in replacing
the auxiliary function $w=a^{(1)}-a^{(2)}$ with $w=a^{(1)}-a^{(2)}-r$, where $r\in
{\mathcal R}_2$ is the minimizer of problem \eqref{eq:4.5.1}, and
noticing that $\widehat{\nabla} r=0$.
\end{proof}

\begin{proof}[Proof of Theorem \ref{theo:Main}]
It is convenient to introduce the following auxiliary distances:
\begin{equation}
    \label{eq:1.0bis}
   d = d_{\cal H}( \overline{\Omega \setminus D_1}, \overline{\Omega \setminus
    D_2}),
\end{equation}
\begin{equation}
    \label{eq:1.0ter}
    d_m = \max \left \{ \max_{x \in \partial D_1} \hbox{dist}(x, \overline{\Omega \setminus
    D_2}),  \max_{x \in \partial D_2} \hbox{dist}(x, \overline{\Omega \setminus
    D_1}) \right \}.
\end{equation}
Let $\eta >0$ such that
\begin{equation}
    \label{eq:1.1}
    \max_{i=1,2} \int_{ (\Omega \setminus \overline{\mathcal G}) \setminus
    \overline{D_i}} |\widehat{\nabla} a^{(i)}|^2 \leq \eta.
\end{equation}

\textit{Step $1$. Let us assume $\eta \leq r_0^2 \| \widehat{N}
  \|_{H^{-1/2}(\partial \Omega, \R^2)}^2$. We have
\begin{equation}
    \label{eq:1.2}
    d_m
    \leq
    C r_0
    \left (
    \frac{ \eta   }{ r_0^2 \| \widehat{N}
  \|_{H^{-1/2}(\partial \Omega, \R^2)}^2  }
    \right )^{\frac{1}{\tau}},
\end{equation}
where $\tau$ has been introduced in Proposition \ref{prop:FVRB}
and $C$ is a positive constant only depending on the a priori
data.}

\begin{proof} Without loss of generality, let
$x_0 \in \partial D_1$ such that
\begin{equation}
    \label{eq:2.1}
    \hbox{dist}(x_0, \overline{\Omega \setminus D_2}) =d_m>0.
\end{equation}
Since $B_{d_m}(x_0) \subset D_2 \subset \Omega \setminus
\overline{\mathcal G}$, we have
\begin{equation}
    \label{eq:2.2}
    B_{d_m}(x_0) \cap (\Omega \setminus \overline{D_1}) \subset
    ( \Omega \setminus
\overline{\mathcal G} ) \setminus \overline{D_1}
\end{equation}
and then, by \eqref{eq:1.1},
\begin{equation}
    \label{eq:2.3}
    \int_{  B_{d_m}(x_0) \cap (\Omega \setminus \overline{D_1}) } |\widehat{\nabla} a^{(1)}|^2 \leq \eta.
\end{equation}
Let us distinguish two cases. First, let
\begin{equation}
    \label{eq:3.1}
    d_m < \overline{c}_0 r_0,
\end{equation}
where $\overline{c}_0$ is the positive constant appearing in
Proposition \ref{prop:FVRB}. By applying this proposition, we have

\begin{equation}
    \label{eq:4.1}
    \eta \geq C \left(\frac{d_m}{r_0}\right)^\tau \int_{ B_{r_0}(x_0) \cap (\Omega \setminus \overline{D_1}) } |\widehat{\nabla} a^{(1)}|^2,
\end{equation}
where $C>0$ is a positive constant only depending on $\alpha_0$,
$\gamma_0$, $\Lambda_0$, $\alpha$, $M_0$, $M_1$ and $F$.

By Proposition \ref{prop:LPS}, we have

\begin{equation}
    \label{eq:4.3}
    \eta \geq C\left(\frac{d_m}{r_0}\right)^\tau r_0^2 \|\widehat{N}\|_{H^{-1/2}(\partial\Omega, \R^2)}^2,
\end{equation}
where $C>0$ is a positive constant only depending on $\alpha_0$,
$\gamma_0$, $\Lambda_0$, $\alpha$, $M_0$, $M_1$, $F$, {}from which
we can estimate $d_m$, obtaining \eqref{eq:1.2}.

As second case, let
\begin{equation}
    \label{eq:5.2}
    d_m
    \geq
    \overline{c}_0 r_0.
\end{equation}
By starting again {}from \eqref{eq:2.3}, applying Proposition
\ref{prop:LPS} and recalling $d_m\leq M_1r_0$, we have
\begin{equation}
    \label{eq:5.4}
    d_m
    \leq
    C r_0
    \left (
    \frac{ \eta   }{ r_0^2 \|\widehat{N}\|_{H^{-1/2}(\partial\Omega, \R^2)}^2 }
    \right ),
\end{equation}
where $C>0$ is a positive constant only depending on $\alpha_0$,
$\gamma_0$, $\Lambda_0$, $M_0$, $M_1$, $F$. Since we have assumed
$\eta \leq r_0^2 \| \widehat{N} \|_{H^{-1/2}(\partial \Omega,
\R^2)}^2$, also in this case we obtain \eqref{eq:1.2}.
\end{proof}
\textit{Step $2$. Let us assume $\eta \leq r_0^2 \| \widehat{N}
\|_{H^{-1/2}(\partial \Omega, \R^2)}^2$. We have
\begin{equation}
    \label{eq:9bis.1}
    d
    \leq
    C r_0
    \left (
    \frac{ \eta   }{ r_0^2 \|\widehat{N}\|_{H^{-1/2}(\partial\Omega, \R^2)}^2  }
    \right )^{\frac{1}{\tau_1}},
\end{equation}
with $\tau_1=\max\{\tau,\tau_0\}$ and $C>0$ only depends on
$\alpha_0$, $\gamma_0$, $\Lambda_0$, $\alpha$, $M_0$, $M_1$ and
$F$.}
\begin{proof} We may assume that $d>0$ and there exists $y_0\in \overline{\Omega\setminus
D_1}$ such that
\begin{equation}
    \label{eq:7.1}
\hbox{dist}(y_0,\overline{\Omega\setminus D_2})=d.
\end{equation}
Since $d>0$, we have $y_0\in D_2\setminus D_1$. Let
\begin{equation}
    \label{eq:7.2}
    h=\hbox{dist}(y_0,\partial D_1),
\end{equation}
possibly $h=0$.

There are three cases to consider:

i) $h\leq \frac{d}{2}$;

ii) $h> \frac{d}{2}$, $h\leq \frac{d_0}{2}$;

iii) $h> \frac{d}{2}$, $h> \frac{d_0}{2}$.

\medskip
\noindent Here the number $d_0$, $0<d_0<r_0$, is such that
$\frac{d_0}{r_0}$ only depends on $M_0$, and it is the same
constant appearing in Proposition \ref{prop:4.8.1}. In particular,
Proposition $3.6$ in \cite{A-B-R-V} shows that there exists an
absolute constant $C>0$ such that if $d \leq d_0$, then $d\leq
Cd_m$.

\medskip
\noindent \emph{Case i).}

By definition, there exists $z_0\in \partial D_1$ such that
$|z_0-y_0|=h$. By applying the triangle inequality, we get
$\hbox{dist}\left(z_0, \overline{\Omega\setminus D_2}\right)\geq
\frac{d}{2}$. Since, by definition, $\hbox{dist}\left(z_0,
\overline{\Omega\setminus D_2}\right)\leq d_m$, we obtain $d\leq
2d_m$.

\medskip
\noindent \emph{Case ii).}

It turns out that $d<d_0$ and then, by the above recalled
property, again we have that $d\leq Cd_m$, for an absolute
constant $C$.

\medskip
\noindent \emph{Case iii).}

Let $\widetilde{h}=\min\{h,r_0\}$. We obviously have that
$B_{\widetilde{h}}(y_0)\subset \Omega\setminus \overline{D_1}$ and
$B_d(y_0)\subset D_2$. Let us set
$$
d_1=\min\left\{\frac{d}{2},\frac{\widetilde{c_0} d_0}{4}\right\},
$$
where $\widetilde{c_0}$ is the positive constant appearing in
Proposition \ref{prop:FVRI}. Since $d_1<d$ and
$d_1<\widetilde{h}$, we have that $B_{d_1}(y_0)\subset
D_2\setminus \overline{D_1}$ and therefore $\eta\geq
\int_{B_{d_1}(y_0)}|\widehat{\nabla} a^{(1)}|^2$.

Since $\frac{d_0}{2}<\widetilde{h}$,
$B_{\frac{d_0}{2}}(y_0)\subset \Omega\setminus \overline{D_1}$ so
that we can apply Proposition \ref{prop:FVRI} with $r_1=d_1$,
$\overline{r}=\frac{d_0}{2}$, obtaining $\eta\geq
C\left(\frac{2d_1}{d_0}\right)^{\tau_0}\int_{B_{\frac{d_0}{2}}(y_0)}|\widehat{\nabla}
a^{(1)}|^2$, with $C>0$ only depending on $\alpha_0$, $\gamma_0$,
$\Lambda_0$, $M_0$, $M_1$ and $F$. Next, by Proposition
\ref{prop:LPS}, recalling that $\frac{d_0}{r_0}$ only depends on
$M_0$, we derive that
$$
d_1
    \leq
    C r_0
    \left (
    \frac{ \eta   }{ r_0^2 \| \widehat{M} \|_{H^{-1/2}(\partial \Omega, \R^2)}^2  }
    \right )^{\frac{1}{\tau_0}},
$$
where $C>0$ only depends on $\alpha_0$, $\gamma_0$, $\Lambda_0$,
$M_0$, $M_1$ and $F$. For $\eta$ small enough,
$d_1<\frac{\widetilde{c_0}d_0}{4}$, so that $d_1=\frac{d}{2}$ and
$$
d
    \leq
    C r_0
    \left (
    \frac{ \eta   }{ r_0^2 \| \widehat{M} \|_{H^{-1/2}(\partial \Omega, \R^2)}^2   }
    \right )^{\frac{1}{\tau_0}},
$$
where $C>0$ only depends on $\alpha_0$, $\gamma_0$, $\Lambda_0$,
$M_0$, $M_1$ and $F$. Collecting the three cases, the thesis
follows.
\end{proof}

\textit{Step $3$. We have
\begin{equation}
    \label{eq:step-3}
    d_{\mathcal{H}} (\overline{D_1}, \overline{D_2}) \leq
    \sqrt{1+M_0^2} \ d.
\end{equation}
}
\begin{proof}
The proof is based on purely geometrical arguments, we refer to
\cite[Proof of Theorem 3.1, Step 3]{l:mrv2019}.
\end{proof}

\textit{Conclusion.} By Proposition \ref{prop:4.8.1},
\begin{equation}
    \label{eq:ultima}
    d
    \leq
    C r_0
    \left (\log\left|\log\left(\frac{ \varepsilon   }{  \|\widehat{N} \|_{H^{-1/2}(\partial \Omega, \R^2)}^2  }\right)\right|
    \right )^{-\frac{1}{2\tau_1}},
\end{equation}
with $\tau_1\geq 1$ and $C>0$ only depends on $\alpha_0$,
$\gamma_0$, $\Lambda_0$, $\alpha$, $M_0$, $M_1$ and $F$. By this
first rough estimate, there exists $\varepsilon_0>0$, only depending
on on $\alpha_0$, $\gamma_0$, $\Lambda_0$, $\alpha$, $M_0$, $M_1$
and $F$, such that, if $\varepsilon\leq \varepsilon_0$, then $d\leq
d_0$. Therefore, we can apply the second statement of Proposition
\ref{prop:4.8.1}, obtaining the thesis.
\end{proof}

\section{Generalized Plane Stress problem}
   \label{sec:GPS}

In this section we derive the Generalized Plane Stress (GPS)
problem for the statical equilibrium of a thin elastic plate under
in-plane boundary loads. Our analysis follows the classical
approach of the theory of structures, according to the original
idea introduced by Filon \cite{Filon}. Alternative, more formal
derivations have been proposed to justify the GPS problem. The
interested reader can refer, among others, to the contributions
\cite{Ciarlet-Destuynder-1979}, \cite{l:abp} and \cite{Paroni}.

Let ${\mathcal U}$ be a bounded domain in $\R^2$, and consider the
cylinder $\mathcal{C} = {\mathcal U} \times \left ( - \frac{h}{2},
\frac{h}{2} \right )$ with middle plane ${\mathcal U} \times
\{x_3=0\}$ (which we will simply denote by ${\mathcal U}$ in what
follows) and thickness $h$. Here, $\{O,x_1,x_2,x_3\}$ is a
Cartesian coordinate system, with origin $O$ belonging to the
plane $x_3=0$ and axis $x_3$ orthogonal to ${\mathcal U}$. Such
cylinder is called \textit{plate} if $h$ is small with respect to
the linear dimensions of ${\mathcal U}$, e.g., $h <<
diam({\mathcal U})$.

Let us suppose that the faces ${\mathcal U} \times \{x_3 = \pm \frac{h}{2}\}$
of the plate are free of applied loads, and all external surface
forces acting on the lateral surface $\partial {\mathcal U} \times \left ( -
\frac{h}{2}, \frac{h}{2} \right ) $ lie in planes parallel to the
middle plane ${\mathcal U}$, and are independent of $x_3$. We shall further
assume that body forces vanish in $\mathcal{C}$. The plate is
assumed to be made by linearly elastic isotropic material, with
Lam\'{e} moduli independent of the $x_3$-coordinate, e.g.,
$\lambda=\lambda (x_1,x_2)$, $\mu=\mu(x_1,x_2)$ for every
$(x_1,x_2, 0) \in \overline{{\mathcal U}}$. Moreover, let $\lambda$, $\mu \in
C^{0,1}(\overline{{\mathcal U}})$ and such that $\mu \geq \alpha_0$, $2\mu
+3\lambda \geq \gamma_0$ in $\overline{{\mathcal U}}$, with $\alpha_0$,
$\gamma_0$ positive constants.

Under the above assumptions, the problem of elastostatics consists
in finding a displacement $u$ solution to
\begin{center}
\( {\displaystyle \left\{
\begin{array}{lr}
     T_{ij,j}=0,                                        & \mathrm{in} \ \mathcal{C},
        \vspace{0.25em}\\
     T_{i3}=0,                                          & \mathrm{on} \ {\mathcal U} \times \{x_3 = \pm\frac{h}{2} \},
          \vspace{0.25em}\\
     T_{\alpha \beta}n_\beta = \widehat{t}_\alpha,      & \mathrm{on} \ \partial {\mathcal U} \times \left ( - \frac{h}{2},
                                                        \frac{h}{2} \right),
        \vspace{0.25em}\\
        T_{3 \beta}n_\beta = 0,                         & \mathrm{on} \ \partial {\mathcal U} \times \left ( - \frac{h}{2},
                                                        \frac{h}{2} \right),
        \vspace{0.25em}\\
         T_{ij} = 2\mu E_{ij} + \lambda (E_{kk})\delta_{ij},             & \mathrm{in} \ \mathcal{C},
        \vspace{0.25em}\\
         E_{ij} = \frac{1}{2} \left (u_{i,j} + u_{j,i} \right  ),                 & \mathrm{in} \ \mathcal{C},
        \vspace{0.25em}\\
\end{array}
\right. } \) \vskip -10.3em
\begin{eqnarray}
& & \label{eq:GPS.2.1}\\
& & \label{eq:GPS.2.2}\\
& & \label{eq:GPS.2.3}\\
& & \label{eq:GPS.2.4}\\
& & \label{eq:GPS.2.5}\\
& & \label{eq:GPS.2.6}
\end{eqnarray}
\end{center}
where the force field $\widehat{t}=(\widehat{t}_1,
\widehat{t}_2,0)$, with
$\widehat{t}_\alpha=\widehat{t}_\alpha(x_1,x_2)$, $\alpha=1,2$,
assigned on $\partial {\mathcal U} \times \left ( -
\frac{h}{2},\frac{h}{2} \right)$ satisfies the compatibility
conditions
\begin{equation}
  \label{eq:GPS.2.7}
  \int_{\partial {\mathcal U} \times \left ( -
\frac{h}{2},\frac{h}{2} \right)} \widehat{t}=0, \qquad
\int_{\partial {\mathcal U} \times \left ( - \frac{h}{2},\frac{h}{2} \right)}
x \times \widehat{t}=0,
\end{equation}
see, for example, \cite[\S $45$]{Gurtin}. The above boundary value
problem is called \textit{plane problem of elastostatics}. It is
known that, under our assumptions and for $\widehat{t}_\alpha \in
H^{ -\frac{1}{2}} (\partial \mathcal U, \R^2)$, $\alpha=1,2$,
there exists a solution $u \in H^1 (\mathcal{C}, \R^3)$ which is
unique up to an infinitesimal rigid displacement $r(x)=a + b
\times x$, with $a, \ b \in \R^3$ constant vectors.

We now formulate the Generalized Plane Stress (GPS) problem
associated to \eqref{eq:GPS.2.1}--\eqref{eq:GPS.2.6}. The GPS
problem is a two-dimensional boundary value problem formulated in
terms of the thickness averages of $u$, $E$ and $T$, under the a
priori assumption
\begin{equation}
  \label{eq:GPS.3.1}
  T_{33}=0,                 \quad \mathrm{in} \ \mathcal{C}.
\end{equation}
For a physically plausible justification of the above assumption
under the hypothesis of small $h$, we refer to \cite[\S
$67$]{Sokolnikoff} and to the paper \cite{Filon} by Filon, who
first derived the GPS problem.

Given a function $f: \mathcal{C} \rightarrow \R^3$, $f \in
H^1(\mathcal{C})$, let us define the function $\widetilde{f}:
\mathcal{C} \rightarrow \R^3$ as follows:
\begin{center}
\( {\displaystyle \left\{
\begin{array}{lr}
     \widetilde{f}_1(x_1,x_2,x_3)=f_1(x_1,x_2,-x_3),
        \vspace{0.25em}\\
     \widetilde{f}_2(x_1,x_2,x_3)=f_2(x_1,x_2,-x_3),
          \vspace{0.25em}\\
     \widetilde{f}_3(x_1,x_2,x_3)=-f_3(x_1,x_2,-x_3).
        \vspace{0.25em}\\
\end{array}
\right. } \) \vskip -6.0em
\begin{eqnarray}
& & \label{eq:symm.1.1}\\
& & \label{eq:symm.1.2}\\
& & \label{eq:symm.1.3}
\end{eqnarray}
\end{center}
By definition of the plane problem, if $u$ is a solution to
\eqref{eq:GPS.2.1}--\eqref{eq:GPS.2.6}, then also $\widetilde{u}$
is a solution of the same problem. Moreover, $(u-\widetilde{u})$
is a solution to \eqref{eq:GPS.2.1}--\eqref{eq:GPS.2.6} with
$\widehat{t}=0$ and, therefore, $ (u-\widetilde{u}) \in
\mathcal{R}_3$. Noticing that
$(u_1-\widetilde{u}_1)|_{x_3=0}=(u_2-\widetilde{u}_2)|_{x_3=0}=0$,
we have $u-\widetilde{u}=a_3 e_3 + (b_1 e_1 + b_2 e_2 )\times
\sum_{i=1}^3 x_i e_i$, with $a_3, b_1, b_2 \in \R$. Now, it is
easy to see that, choosing $r' \in \mathcal{R}_3$ as $r'=
\sum_{i=1}^3 a_i' e_i + \sum_{i=1}^3 b_i' e_i \times \sum_{i=1}^3
x_i e_i$, with $a_3' = - \frac{a_3}{2}$, $b_1'= - \frac{b_1}{2}$,
$b_2' = - \frac{b_2}{2}$, the solution $u+r'$ to
\eqref{eq:GPS.2.1}--\eqref{eq:GPS.2.6} satisfies the condition
$u+r'= \widetilde{(u+r')}$, for every $a_1'$, $a_2'$, $b_3' \in
\R$.

We next introduce the \textit{thickness average} $\overline{f}$
 of a function $f: \mathcal{C}
\rightarrow \R^3$, $\overline{f} : U \rightarrow \R$, defined as
\begin{equation}
  \label{eq:GPS.4.1}
    \overline{f}(x_1,x_2)= \frac{1}{h} \int_{- \frac{h}{2}}^{
    \frac{h}{2}} f(x_1,x_2,x_3)dx_3.
\end{equation}
Taking into account that the thickness average of an $x_3$-odd
function is zero, and the $x_3$-derivative of an $x_3$-even
function is $x_3$-odd, for every point $(x_1,x_2) \in {\mathcal U}$ we have
\begin{center}
\( {\displaystyle \left\{
\begin{array}{lr}
    \overline{u}_3=\overline{E}_{\alpha 3}=\overline{T}_{\alpha
    3}=0, \quad \alpha=1,2,
        \vspace{0.25em}\\
     \overline{E}_{\alpha \beta} = \frac{1}{2} \left (
     \overline{u}_{\alpha,\beta} + \overline{u}_{\beta,\alpha} \right ),
     \quad \alpha, \beta=1,2,
          \vspace{0.25em}\\
     \overline{T}_{\alpha \beta} = 2\mu \overline{E}_{\alpha
     \beta}+ \lambda ( \overline{E}_{\gamma
     \gamma} + \overline{E}_{33} )\delta_{\alpha \beta}, \quad \alpha, \beta=1,2,
        \vspace{0.25em}\\
     \overline{T}_{33} = 2\mu \overline{E}_{33}+ \lambda ( \overline{E}_{\gamma
     \gamma} + \overline{E}_{33} ),
        \vspace{0.25em}\\
\end{array}
\right. } \) \vskip -7.1em
\begin{eqnarray}
& & \label{eq:GPS.4.2}\\
& & \label{eq:GPS.4.3}\\
& & \label{eq:GPS.4.4}\\
& & \label{eq:GPS.4.5}
\end{eqnarray}
\end{center}
where the solution $u+r'$ is denoted by $u$. Using the a priori
assumption \eqref{eq:GPS.3.1} in \eqref{eq:GPS.4.5}, the function
$\overline{E}_{33}$ can be expressed in terms of
$\overline{E}_{\gamma \gamma}$, and the \textit{two-dimensional}
constitutive equation can be written as
\begin{equation}
  \label{eq:GPS.4.6}
    \overline{T}_{\alpha \beta} = 2\mu \overline{E}_{\alpha
     \beta}+ \lambda^* \overline{E}_{\gamma
     \gamma} \delta_{\alpha \beta},
\end{equation}
with
\begin{equation}
  \label{eq:GPS.4.7}
    \lambda^* = \frac{2\mu\lambda}{\lambda + 2\mu}.
\end{equation}
Integrating on the thickness in
\eqref{eq:GPS.2.1}--\eqref{eq:GPS.2.6}, and neglecting those
equations which yield to identities, we obtain the averaged
equations of equilibrium and the corresponding boundary
conditions, and $\overline{u} \in H^1({\mathcal U}, \R^2)$ is a solution to
\begin{center}
\( {\displaystyle \left\{
\begin{array}{lr}
     \overline{T}_{\alpha\beta,\beta}=0,                                        & \mathrm{in} \ {\mathcal U},
        \vspace{0.25em}\\
     \overline{T}_{\alpha \beta}n_\beta = \widehat{t}_\alpha,      & \mathrm{on} \ \partial {\mathcal U},
        \vspace{0.25em}\\
     \overline{T}_{\alpha\beta} = 2\mu \overline{E}_{\alpha\beta} +  \lambda^* (\overline{E}_{\gamma\gamma})\delta_{\alpha\beta},             & \mathrm{in} \ \mathcal U,
        \vspace{0.25em}\\
     \overline{E}_{\alpha\beta} = \frac{1}{2} \left (\overline{u}_{\alpha,\beta} + \overline{u}_{\beta,\alpha} \right  ),           & \mathrm{in} \
     {\mathcal U},
        \vspace{0.25em}\\
\end{array}
\right. } \) \vskip -7.3em
\begin{eqnarray}
& & \label{eq:GPS.5.1}\\
& & \label{eq:GPS.5.2}\\
& & \label{eq:GPS.5.3}\\
& & \label{eq:GPS.5.4}
\end{eqnarray}
\end{center}
where the force field $\widehat{t}=\widehat{t}_1 e_1 +
\widehat{t}_2 e_2$ applied on $\partial {\mathcal U}$ satisfies the
compatibility conditions
\begin{equation}
  \label{eq:GPS.comp}
   \int_{\partial {\mathcal U}} \widehat{t}=0, \quad \int_{\partial {\mathcal U}} x \times
   \widehat{t}=0.
\end{equation}
Let us notice that the constitutive equation \eqref{eq:GPS.5.3}
can be written as
\begin{equation}
  \label{eq:GPS.5.5}
    \overline{T}_{\alpha\beta} =
    \frac{E}{1-\nu^2}\left((1-\nu)\overline{E}_{\alpha\beta}+\nu(\overline{E}_{\gamma\gamma})\delta_{\alpha\beta}\right),
\end{equation}
with
\begin{equation}
  \label{eq:GPS.5.6}
    \mu = \frac{E}{2(1+\nu)}, \qquad \lambda = \frac{\nu
    E}{(1+\nu)(1-2\nu)},
\end{equation}
where $E$, $\nu$ are the Young's modulus and the Poisson's
coefficient of the material, respectively. Finally, by defining
\begin{equation}
  \label{eq:GPS.6.1}
    a_\alpha = \overline{u}_\alpha,
    \quad
    \epsilon_{\alpha\beta}= \overline{E}_{\alpha\beta}= \widehat{\nabla}a,
    \quad
    N_{\alpha \beta}=h
    \overline{T}_{\alpha\beta},
    \quad
    \widehat{N}_\alpha = h\widehat{t}_\alpha,
    \quad
    \alpha, \beta =1,2,
\end{equation}
we obtain the GPS problem
\begin{center}
\( {\displaystyle \left\{
\begin{array}{lr}
     N_{\alpha\beta,\beta}=0,
      & \mathrm{in}\ {\mathcal U},
        \vspace{0.25em}\\
      N_{\alpha\beta}n_\beta=\widehat{N}_\alpha, & \mathrm{on}\ \partial {\mathcal U},
          \vspace{0.25em}\\
      N_{\alpha\beta}=\frac{Eh}{1-\nu^2}\left((1-\nu)\epsilon_{\alpha\beta}+\nu(\epsilon_{\gamma\gamma})\delta_{\alpha\beta}\right), &\mathrm{in}\ \overline{{\mathcal U}},
          \vspace{0.25em}\\
        \epsilon_{\alpha\beta}=\frac{1}{2}(a_{\alpha,\beta}+a_{\beta,\alpha}), &\mathrm{in}\
        {\mathcal U},
          \vspace{0.25em}\\
\end{array}
\right. } \) \vskip -7.5em
\begin{eqnarray}
& & \label{GPS.6.7.1}\\
& & \label{GPS.6.7.2}\\
& & \label{GPS.6.7.3}\\
& & \label{GPS.6.7.4}
\end{eqnarray}
\end{center}
with
\begin{equation}
  \label{GPS.6.7.5}
  \int_{\partial {\mathcal U}}\widehat{N}\cdot r =0, \quad \hbox{for every } r\in {\mathcal
  R_2}.
\end{equation}

\section*{Acknowledgements}
The authors wish to thank Jenn-Nan Wang for fruitful 
discussions on the subject of this work. Antonino Morassi is supported by PRIN 2015TTJN95 ``Identificazione
e diagnostica di sistemi strutturali complessi''. Edi Rosset and
Sergio Vessella are supported by Progetti GNAMPA 2017 and 2018,
Istituto Nazionale di Alta Matematica (INdAM). Edi Rosset is
supported by FRA $2016$ ``Problemi inversi, dalla stabilit\`a alla
ricostruzione", Universit\`a degli Studi di Trieste.

\end{document}